\documentclass{amsart}
\usepackage{graphics}

\newtheorem{theorem}{Theorem}[section]
\newtheorem{lemma}[theorem]{Lemma}
\newtheorem{cor}[theorem]{Corollary}

\theoremstyle{remark}

\numberwithin{equation}{section}

\begin{document}
\newcommand{\beqs}{\begin{equation*}}
\newcommand{\eeqs}{\end{equation*}}
\newcommand{\beq}{\begin{equation}}
\newcommand{\eeq}{\end{equation}}
\newcommand{\bin}[2]{\genfrac{(}{)}{0pt}{0}{#1}{#2}}
\newcommand\mymod[1]{\mbox{ mod}\ {#1}}
\newcommand\smymod[1]{\mbox{\scriptsize mod}\ {#1}}
\newcommand\mylabel[1]{\label{#1}}
\newcommand\eqn[1]{(\ref{eq:#1})}

\title[\tiny Ternary Quadratic Forms, Modular Equations, Certain Positivity Conjectures]
{Ternary Quadratic Forms, Modular Equations \\ and Certain Positivity Conjectures.}

\author{Alexander Berkovich}
\address{Department of Mathematics, University of Florida, Gainesville,
Florida 32611-8105}
\email{alexb@ufl.edu}
\thanks{Research of the first author was supported in part by NSA grant
H98230-09-1-0051}

\author{William C. Jagy}
\address{Math.Sci.Res.Inst.,1000 Centenial Drive,Berkeley, CA94720}
\email{jagy@msri.org}


\date{June 16th, 2009}

\dedicatory{Dedicated to the memory of Professor Alladi Ramakrishnan.}


\begin{abstract}

We show that many of Ramanujan's modular equations of degree $3$ can be interpreted in terms of integral 
ternary quadratic forms. This way we establish that for any $n\in\mathbf N$
\beqs
\begin{split}
& \left|\left\{n=\frac{x(x+1)}{2}+y^2+z^2:x,y,z\in\mathbf Z\right\}\right|\geq \\
& \left|\left\{n=\frac{x(x+1)}{2}+3y^2+3z^2:x,y,z\in\mathbf Z\right\}\right|,
\end{split}
\eeqs
just to name one among many similar ``positivity'' results of this type.
In particular, we prove the recent conjecture of H. Yesilyurt and the first author, stating that for any $n\in\mathbf N$
\beqs
\begin{split}
& \left|\left\{n=\frac{x(x+1)}{2}+y^2+z^2:x,y,z\in\mathbf Z\right\}\right|\geq \\
& \left|\left\{n=\frac{x(x+1)}{2}+7y^2+7z^2:x,y,z\in\mathbf Z\right\}\right|.
\end{split}
\eeqs
We prove a number of identities for certain ternary forms with discriminant $144,400,784,3600$
by converting every ternary identity into an identity for the appropriate $\eta$-quotients.
In the process we discover and prove a few new modular equations of degree $5$ and $7$.
For any square free odd integer $S$ with prime factorization $p_1\ldots p_r$, we define the $S$-genus 
as a union of $2^r$ specially selected genera of ternary quadratic forms, all with discriminant $16S^2$.
This notion of $S$-genus arises naturally in the course of our investigation. It entails an interesting 
injection from genera of binary quadratic forms with discriminant $-8S$ to genera of ternary quadratic 
forms with discriminant $16S^2$.

\end{abstract}

\maketitle

\section{Introduction} \label{sec:intro}
\medskip

Alladi Ramakrishnan's visits to Gainesville were always memorable. His
interests were diverse and his passion for science was truly amazing.
He was a very open man, always happy to make new friends. He had so many
stories to tell. His family was a true pillar of strength for him  and
in turn he was devoted to them.
\medskip

Ramanujan's general theta-function $f(a,b)$ is defined by
\beq
f(a,b)=\sum_{n=-\infty}^{\infty}a^{\frac{(n-1)n}{2}}b^{\frac{(n+1)n}{2}},\quad |ab|<1.
\mylabel{eq:1.1}
\eeq
In Ramanujan's notation, the celebrated Jacobi triple product identity takes the shape
\beq
f(a,b)=(-a;ab)_\infty(-b;ab)_\infty(ab;ab)_\infty,\quad |ab|<1,
\mylabel{eq:1.2}
\eeq
with 
\beqs
(a;q)_\infty:=\prod_{j\ge 0}(1-aq^j).
\eeqs
It is always assumed that $|q|<1$. The following four special cases will play a prominent role in our narrative
\beq
\phi(q):=f(q,q)=\sum_{n=-\infty}^{\infty}q^{n^2},
\mylabel{eq:1.3}
\eeq
\beq
\psi(q):=f(q,q^3)=\sum_{n\ge 0}q^{\frac{(n+1)n}{2}},
\mylabel{eq:1.4}
\eeq
\beq
f(q,q^2)=\sum_{n=-\infty}^{\infty}q^{\frac{(3n+1)n}{2}},
\mylabel{eq:1.5}
\eeq
\beq
f(q,q^5)=\sum_{n=-\infty}^{\infty}q^{(3n+2)n}.
\mylabel{eq:1.6}
\eeq
Using \eqn{1.2} it is not hard to derive the product representation formulas
\beq
\phi(q)=\frac{E(q^2)^5}{E(q)^2E(q^4)^2},
\mylabel{eq:1.7}
\eeq
\beq
\psi(q)=\frac{E(q^2)^2}{E(q)},
\mylabel{eq:1.8}
\eeq
\beq
f(q,q^2)=\frac{E(q^3)^2E(q^2)}{E(q^6)E(q)},
\mylabel{eq:1.9}
\eeq
\beq
f(q,q^5)=\frac{E(q^{12})E(q^3)E(q^2)^2}{E(q^6)E(q^4)E(q)},
\mylabel{eq:1.10}
\eeq
where
\beqs
E(q):=\prod_{j\ge 1}(1-q^j).
\eeqs
Combining \eqn{1.7} and \eqn{1.8} we see that
\beq
\psi(q)^2=\phi(q)\psi(q^2). 
\mylabel{eq:1.11}
\eeq
Note that the Dedekind  $\eta(z)$ is related to $E(q)$ as
\beq
\eta(z)=q^{\frac{1}{24}}E(q),\quad\mbox{if\;\;}q=exp(2\pi iz)\quad\mbox{with\;\;}Im(z)>0.
\mylabel{eq:1.12}
\eeq
And so $\phi(q)$, $q^{\frac{1}{8}}\psi(q)$, $q^{\frac{1}{24}}f(q,q^2)$, $q^{\frac{1}{3}}f(q,q^5)$
all have $\eta$-quotient representations.

Following \cite{BeY}, we say that a $q$-series is positive if its power series coefficients are all non-negative.
We define $P[q]$ to be the set of all such series. It is plain that $\phi(q)$, $\psi(q)$, $f(q,q^2)$, $f(q,q^5)$, 
$\frac{1}{E(q)}$ (and their products) are in $P[q]$.
However, it is not at all obvious that
\beq
\psi(q)(\phi(q)^2-\phi(q^7)^2)\in P[q].
\mylabel{eq:1.13}
\eeq
Motivated by their studies of $7$-core partitions, H. Yesilyurt and the first author conjectured \eqn{1.13} in ((6.2), \cite{BeY}).
The reader should be cautioned that other similar conjectures there: (6.1), (6.3) and (6.4) are false.
What makes \eqn{1.13} somewhat  non-trivial is the fact that it is not true that $\phi(q)^2-\phi(q^7)^2\in P[q]$.
However, in Section 3 we shall prove
\begin{theorem}\label{t1}
The following identities are true
\begin{align}
\psi(q)(\phi(q)^2-\phi(q^7)^2) & = 4q\psi(q^2)\psi(q^7)\phi(q^7)+8q^2\psi(q^{14})\psi(q^3)\phi(q^{21}) \nonumber \\
&+8q^4\psi(q^{14})f(q,q^2)f(q^7,q^{35}),\mylabel{eq:1.14} \\
\intertext{and}
7\phi(q^7)^2\psi(q^7) & = 8q^2\psi(q^2)\psi(q^{21})\phi(q^3)+8\psi(q^2)f(q,q^5)f(q^7,q^{14}) \nonumber \\
&-4q\psi(q^{14})\psi(q)\phi(q)- \psi(q^7)\phi(q)^2.
\mylabel{eq:1.15}
\end{align}
\end{theorem}
\noindent
One does not have to be very perceptive to deduce that the right hand side (RHS) of \eqn{1.14} is in $P[q]$.
Hence, \eqn{1.13} follows. Our proof of \eqn{1.14} makes naive use of the theory of modular forms.
We employ certain of Ramanujan's modular equations of degree $7$ \cite{Bern} to derive \eqn{1.15} from \eqn{1.14}.
Also in Section 3 we provide the following beautiful interpretation of Theorem 1.1 in terms of integral ternary quadratic forms
\begin{theorem}\label{t2}
If $M\equiv 1\mymod 8$, then
\beq
\begin{split}
(1,8,8,0,0,0)(M) & = (1,14,14,0,0,0)(M) \\
& + 2(2,7,14,0,0,0)(M)+4(3,5,14,0,0,2)(M).
\mylabel{eq:1.16}
\end{split}
\eeq
If $M\equiv 1\mymod 8$ and $7| M$, then
\beq
\begin{split}
7(1,8,8,0,0,0)\left(\frac{M}{7^2}\right) = & -(1,14,14,0,0,0)(M) \\
& - 2(2,7,14,0,0,0)(M)+4(3,5,14,0,0,2)(M),
\mylabel{eq:1.17}
\end{split}
\eeq
where here and everywhere
\beqs
(a,b,c,d,e,f)(M):=|\{M=ax^2+by^2+cz^2+dyz+ezx+fxy:x,y,z\in\mathbf Z\}|.
\eeqs
\end{theorem}
\noindent
We have the following
\begin{cor}\label{c1}
If $M\equiv 1\mymod 8$ and $(M|7)=1$, then
\beqs
(1,8,8,0,0,0)(M)=(1,14,14,0,0,0)(M)+2(2,7,14,0,0,0)(M).
\eeqs
If $M\equiv 1\mymod 8$, $7||M$, then
\beqs
(1,8,8,0,0,0)(M)=2(1,14,14,0,0,0)(M)+4(2,7,14,0,0,0)(M).
\eeqs
\end{cor}
\noindent
A few remarks are in order.
We use the convention that Jacobi's symbol $(M|a)=0$, whenever $(M,a)>1$.
The notation $p||n$ means that $p|n$ but it is not true that $p^2|n$.
A slightly different version of this corollary was communicated to us by Benjamin Kane. His observation
was crucial to our investigation. We understand that Kane used Siegel's weighted average theorem \cite{Iwan}
together with some local calculations of Jones \cite{Jones}.
We note that the Corollary 1.3 has a twin
\begin{cor}\label{c2}
If $M\equiv 1\mymod 8$ and $(M|7)=-1$, then
\beqs
(1,8,8,0,0,0)(M)= 4(3,5,14,0,0,2)(M).
\eeqs
If $M\equiv 1\mymod 8$, $7||M$, then
\beqs
(1,8,8,0,0,0)(M)=8(3,5,14,0,0,2)(M).
\eeqs
\end{cor}
There is nothing very special about the exponent $7$ in \eqn{1.13}. In fact, we plan to prove that for any $S\in\mathbf N$
\beq
\psi(q)(\phi(q)^2-\phi(q^S)^2)\in P[q].
\mylabel{eq:1.18}
\eeq

In this paper we discuss in great detail $S=3,5,7,15$.
In these cases we will construct and prove $\eta$-quotient identities that imply appropriate positivity results.
For $S=3$ and $5$ our identities can be written concisely as modular equations of degree $3$ and $5$, respectively.
Ramanujan found an astounding number of modular equations of degree $3$ and $5$. These are collected and proven  in   \cite{Bern}. The results there are sufficient to prove everything needed for our treatment of $S=3$, $5$ 
in Section 2 and Section 3. In Section 4 we prove our Theorem 1.1 and Theorem 1.2. Section 5 deals with the $S=15$ case. 
In Section 6 we define an injective map from genera of binary quadratic forms to genera of ternary
quadratic forms. This map allows us to introduce a very useful notion of $S$-genus.

\bigskip
\section{Ramanujan's modular equations of degree $3$ and associated identities for ternary quadratic forms with discriminant  $144$.} \label{sec:144}
\medskip

Following Ramanujan, we define the multiplier $m$ of degree $n$ as
\beq
m:=m(n,q)=\frac{\phi(q)^2}{\phi(q^n)^2},
\mylabel{eq:2.1}
\eeq
and
\beq
\alpha:=\alpha(q)=1-\frac{\phi(-q)^4}{\phi(q)^4},
\mylabel{eq:2.2}
\eeq
\beq
\beta:=\beta(n,q)=\alpha(q^n).
\mylabel{eq:2.3}
\eeq
We often say that $\beta$ has degree $n$ over $\alpha$. We also call an algebraic relation connecting
$m$, $\alpha$ and $\beta$ a modular equation of degree $n$.
It is well known, page 40, \cite{Bern} that
\beq
\phi(q)=\phi(q^4)+2q\psi(q^8), 
\mylabel{eq:2.4}
\eeq
\beq
\phi(q)^4-\phi(-q)^4=16q\psi(q^2)^4.
\mylabel{eq:2.5}
\eeq
The last equation implies another formula for $\alpha$
\beq
\alpha=16q\frac{\psi(q^2)^4}{\phi(q)^4},  
\mylabel{eq:2.6}
\eeq
which will come in handy later.
Pages 230---237 in \cite{Bern} contain an impressive collection of $15$ of Ramanujan's modular equations of degree $3$ together with succinct proofs. In particular, one can find there what amounts to the following
\begin{lemma}\label{l1}
If
\beqs
\alpha=\frac{p(2+p)^3}{(1+2p)^3},
\eeqs
then 
\beqs
\beta(3,q)=\frac{p^3(2+p)}{(1+2p)},
\eeqs
and
\beqs 
m(3,q)=1+2p.
\eeqs
\end{lemma}
We comment that this lemma is a very efficient tool for verifying any modular equations of degree $3$. 
In particular, we see that
\beq
m-1=2\frac{\beta^{\frac{3}{8}}}{\alpha^{\frac{1}{8}}}.
\mylabel{eq:2.7}
\eeq
From \eqn{1.11}, \eqn{2.1}, \eqn{2.3}, \eqn{2.6} and \eqn{2.7}, it is readily shown that
\beq
\frac{\phi(q)^2}{\phi(q^3)^2}-1=4q\frac{\psi(q)\psi(q^3)\psi(q^6)}{\psi(q^2)\phi(q^3)^2}.
\mylabel{eq:2.8}
\eeq
Clearly, this theta-function identity can be rewritten as
\beq
\psi(q^2)\phi(q)^2=\psi(q^2)\phi(q^3)^2+4q\psi(q)\psi(q^3)\psi(q^6).
\mylabel{eq:2.9}
\eeq
Next, we multiply both sides of \eqn{2.9} by $\frac{\psi(q)}{\psi(q^2)}$ and employ \eqn{1.11} again to deduce that
\beq
\psi(q)\phi(q)^2 =\psi(q)\phi(q^3)^2+4q\psi(q^3)\psi(q^6)\phi(q).
\mylabel{eq:2.10}
\eeq
The truth of
\beq
\psi(q)(\phi(q)^2-\phi(q^3)^2)\in P[q],
\mylabel{eq:2.11}
\eeq
and of
\beq
\psi(q^2)(\phi(q)^2-\phi(q^3)^2)\in P[q]
\mylabel{eq:2.12}
\eeq
is now evident.
Remarkably, we can interpret \eqn{2.9} and \eqn{2.10} in terms of integral ternary quadratic forms with discriminant $144$.
We remind the reader that the discriminant of a ternary form $ax^2+by^2+cz^2+dyz+ezx+fxy$ is defined as 
\beqs
\frac{1}{2}\det
\begin{bmatrix} 2a & f & e \\ f & 2b & d \\ e & d & 2c
\end{bmatrix}.
\eeqs 
To this end we define a sifting operator $S_{t,s}$ by its action on power series as follows
\beq
S_{t,s} \sum_{n\ge0}c(n)q^n=\sum_{k\ge 0}c(tk+s)q^k,
\mylabel{eq:2.13}
\eeq
where $t$, $s$ are integers and $0\leq s < t$. Observe that \eqn{2.4} implies that
\beq
S_{8,1} \phi (q)\phi(q^8)^2=2\psi(q)\phi(q)^2, 
\mylabel{eq:2.14}
\eeq
\beq
S_{8,1}\phi(q)\phi(q^6)^2=2\psi(q)\phi(q^3)^2,
\mylabel{eq:2.15}
\eeq
\beq
S_{8,1}\phi (q^2)\phi(q^3)\phi(q^6) = 4q\phi(q) \psi(q^3) \psi(q^6).
\mylabel{eq:2.16}
\eeq
Employing \eqn{2.14}, \eqn{2.15}, \eqn{2.16} together with \eqn{2.10} we see that
\beq
S_{8,1}\left(\phi(q)\phi(q^8)^2-\phi(q)\phi(q^6)^2-2\phi(q^2)\phi(q^3)\phi(q^6)\right)=0.
\mylabel{eq:2.17}
\eeq
But the above is nothing else but the statement that
\beq
(1,8,8,0,0,0)(M)=(1,6,6,0,0,0)(M)+2(2,3,6,0,0,0)(M),
\mylabel{eq:2.18}
\eeq
for any $M\equiv 1\mymod 8$.
Actually, with the aid of \eqn{2.4} we can easily check that
\beq
S_{8,r}\phi(q)\phi(q^8)^2=\frac{1}{3}S_{8,r}\phi(q)^3 
\mylabel{eq:2.19}
\eeq
with $r=1,7$. 
Hence, \eqn{2.18} may be stated as
\beq
\frac{1}{3}(1,1,1,0,0,0)(M)=(1,6,6,0,0,0)(M)+2(2,3,6,0,0,0)(M),
\mylabel{eq:2.20}
\eeq
for any $M\equiv 1\mymod 8$. 
It is very likely that modular equation \eqn{2.7} was known to Legendre and Jacobi. 
Surprisingly, the quadratic form interpretation given in \eqn{2.20} above appears to be new.
We note that the two ternary forms $x^2+6y^2+6z^2$ and $2x^2+3y^2+6z^2$ on
the right of \eqn{2.20} have the same discriminant $=144$.
Moreover, these two forms have class number $=1$. This means that these
forms belong to different genera and that they are both regular \cite{Dick}, \cite{Jones}, \cite{JKS}.
Moreover, it is easy to see that $(-n_1|3)=-1$ for any integer $n_1$ represented by $x^2+6y^2+6z^2$, $\gcd(n_1,3)=1$ 
and that $(-n_2|3 )=1$ for any integer $n_2$ represented by $2x^2+3y^2+6z^2$, $\gcd(n_2,3)=1$.
We remark that the appearance of at least two genera with the same discriminant is the salient feature of all our ternary form identities. Somewhat anticipating developments in Section 6 we would like to comment that one can obtain the two
ternary forms on the right of \eqn{2.20} starting with binary forms of discriminant $-24$. There are just two genera of
binary quadratic forms with this discriminant: the (equivalence) class of $x^2+6y^2$ and the class of $2x^2+3y^2$ 
(See \cite{Cox}, pages 52--54). All we need to do to obtain our desired ternaries is to add $6z^2$ to $x^2+6y^2$ and $2x^2+3y^2$, respectively. Actually, we can extend \eqn{2.20} a bit as
\beq
\frac{1}{3}(1,1,1,0,0,0)(M)=(1,6,6,0,0,0)(M)+2(2,3,6,0,0,0)(M),
\mylabel{eq:2.21}
\eeq
for any $M\equiv 1,2\mymod 4$. To this end we make repeated use of \eqn{2.4} and confirm that
\beq
S_{4,1}\phi(q)^3=6\psi(q^2)\phi(q)^2, 
\mylabel{eq:2.22}
\eeq
\beq
S_{4,1}\phi(q)\phi(q^6)^2=2\psi(q^2)\phi(q^3)^2, 
\mylabel{eq:2.23}
\eeq
\beq
S_{4,1}\phi(q^2)\phi(q^3)\phi(q^6)=4q\psi(q)\psi(q^3)\psi(q^6) 
\mylabel{eq:2.24}
\eeq
\beq
S_{4,2}\phi(q)^3 = 12\phi(q)\psi(q^2)^2, 
\mylabel{eq:2.25}
\eeq
\beq
S_{4,2}\phi(q)\phi(q^6)^2=4q\psi(q^6)^2\phi(q),
\mylabel{eq:2.26}
\eeq
\beq
S_{4,2}\phi(q^2)\phi(q^3)\phi(q^6)=2\psi(q)\psi(q^3)\phi(q^3).
\mylabel{eq:2.27}
\eeq
Next, we combine \eqn{2.9} and \eqn{2.22}--\eqn{2.24}, to arrive at
\beqs
S_{4,1}\left(\frac{1}{3}\phi(q)^3-\phi(q)\phi(q^6)^2-2\phi(q^2)\phi(q^3)\phi(q^6)\right)=0,
\eeqs 
which is, essentially, the $M\equiv 1\mymod 4$ case in \eqn{2.21}. To see that \eqn{2.21} is also valid when 
$M\equiv 2\mymod 4$ we again use Lemma 2.1 to verify our next modular equation of degree $3$
\beq
m =\frac{\beta^{\frac{1}{2}}}{\alpha^{\frac{1}{2}}}+2\frac{\beta^{\frac{1}{8}}}{\alpha^{\frac{3}{8}}}.
\mylabel{eq:2.28}
\eeq
Indeed, expressing everything in terms of $p$ and simplifying, we obtain the trivial identity
\beqs
1+2p=\frac{p(1+2p)}{2+p}+ 2\frac{1+2p}{2+p}.
\eeqs
Hence, the proof of \eqn{2.28} is complete.
The theta-function identity associated with \eqn{2.28} takes the pleasant form
\beq
\phi(q)\psi(q^2)^2 =\psi(q)\psi(q^3)\phi(q^3)+ q\phi(q)\psi(q^6)^2.
\mylabel{eq:2.29}
\eeq
Hence,
\beqs
\phi(q)(\psi(q^2)^2- q\psi(q^6)^2)\in P[q].
\eeqs
We observe that that $\psi(q^2)^2-q\psi(q^6)^2 \not\in P[q]$.
Again, we combine \eqn{2.25}--\eqn{2.27} and \eqn{2.29} to obtain 
\beqs
S_{4,2}\left(\frac{1}{3}\phi(q)^3-\phi(q)\phi(q^6)^2-2\phi(q^2)\phi(q^3)\phi(q^6)\right)=0,
\eeqs
which is essentially the $M\equiv 2\mymod 4$ case in \eqn{2.21}.
This is not the end of the story, however. We discovered that \eqn{2.18} has an attractive companion
\beq
3(1,8,8,0,0,0)\left(\frac{M}{3^2}\right)=-(1,6,6,0,0,0)(M)+2(2,3,6,0,0,0)(M),
\mylabel{eq:2.30}
\eeq
$M\equiv 1\mymod 8$, $3|M$. To prove it we begin with the modular equation
\beq
\frac{3}{m}+1=2\frac{\alpha^{\frac{3}{8}}}{\beta^{\frac{1}{8}}}, 
\mylabel{eq:2.31}
\eeq
which can be routinely verified with the aid of Lemma 2.1.
Next, we use \eqn{2.1}, \eqn{2.3} and \eqn{2.6} to convert \eqn{2.31} into the theta-function identity
\beq
-\psi(q^3)\phi(q)^2+4\psi(q)\psi(q^2)\phi(q^3)=3\psi(q^3)\phi(q^3)^2.
\mylabel{eq:2.32}
\eeq
With a bit of labor we can show that \eqn{2.32} is equivalent to
\beqs
S_{24,9}(-\phi(q)\phi(q^6)^2 +2\phi(q^2)\phi(q^3)\phi(q^6))=6\psi(q^3)\phi(q^3)^2.
\eeqs
Hence,
\beqs
\begin{split}
& -(1,6,6,0,0,0)(24n+9)+2(2,3,6,0,0,0)(24n+9)= \\
& 3\left|\left\{n= 3x^2+ 3y^2+3\frac{(1+z)z}{2}:\quad x,y,z \in \mathbf Z\right\}\right|.
\end{split}
\eeqs
Observe that
\beqs
\begin{split}
& 3\left|\left\{n= 3x^2+ 3y^2+3\frac{(1+z)z}{2}:\quad x,y,z \in\mathbf Z\right\}\right|= \\
& 3\left|\left\{1+\frac{8n}{3}= 8x^2+ 8y^2+ z^2:\quad x,y,z \in \mathbf Z\right\}\right|.
\end{split}
\eeqs
And so we have completed the proof of \eqn{2.30}. We note the following interesting corollary
\beq
(1,6,6,0,0,0)(M)= 2(2,3,6,0,0,0)(M),
\mylabel{eq:2.33}
\eeq
when $M\equiv 1\mymod 8$, $3|| M$.
Recalling \eqn{2.18} we see that
\beq
(1,8,8,0,0,0)(M) =2(1,6,6,0,0,0)(M),
\mylabel{eq:2.34}
\eeq
with  $M\equiv 1\mymod 8$, $3|| M$.
Analogously, \eqn{2.21} has its own companion identity
\beq
(1,1,1,0,0,0)\left(\frac{M}{3^2}\right)=-(1,6,6,0,0,0)(M)+2(2,3,6,0,0,0)(M),
\mylabel{eq:2.35}
\eeq
with $M\equiv 1,2\mymod 4$, $3|M$.
Since the argument is pretty similar, we confine ourselves to the following diagram \\
Lemma 2.1 \\
$\Downarrow$ \\
Modular equation:  
\beqs
3+ m\frac{\alpha^{\frac{1}{2}}}{\beta^{\frac{1}{2}}} = 2 m\frac{\alpha^{\frac{1}{8}}}{\beta^{\frac{3}{8}}}. 
\eeqs 
$\Downarrow$ \\
Theta-function identity
\beqs
\psi(q^2)^2\phi(q^3)-\psi(q)\psi(q^3)\phi(q)+3q\psi(q^6)\psi(q^3)^2=0.
\eeqs 
$\Downarrow$ \\
Ternary identity  
\beqs
(1,1,1,0,0,0)\left(\frac{M}{3^2}\right)=-(1,6,6,0,0,0)(M)+2(2,3,6,0,0,0)(M ),
\eeqs
with $M\equiv 2\mymod 4$, $3|M$. \\
Lemma 2.1 \\
$\Downarrow$ \\
Modular equation 
\beqs
3+m = 2m\frac{\alpha^{\frac{3}{8}}}{\beta^{\frac{1}{8}}}. 
\eeqs
$\Downarrow$ \\
Theta-function identity 
\beqs
\psi(q^6)\phi(q)^2-4\psi(q)\psi(q^2)\psi(q^3) +3\psi(q^3)^2\phi(q^3)=0.
\eeqs
$\Downarrow$ \\
Ternary identity
\beqs
(1,1,1,0,0,0)\left(\frac{M}{3^2}\right)  = -(1,6,6,0,0,0)(M) +2(2,3,6,0,0,0)(M)
\eeqs
with  $M\equiv 1\mymod 4$, $3|M$.
We conclude this section by stating \eqn{2.18}, \eqn{2.30} in a way that would
suggest an elegant and straightforward generalization.
Let $|\mbox{Aut(a,b,c,d,e,f)}|$  denote the number of integral automorphs of a
ternary form $ax^2+by^2+cz^2+dyz+ezx+fxy$.
It is easy to check that  
$$\frac{16}{|\mbox{Aut(1,6,6,0,0,0)}|}=1,$$ 
$$\frac{16}{|\mbox{Aut(2,3,6,0,0,0)}|}=2.$$
And so we can rewrite RHS \eqn{2.18} as a weighted average over two genera.  
This way it becomes
\beq
(1,8,8,0,0,0)(M) =\frac{16(1,6,6,0,0,0)(M)}{|\mbox{Aut(1,6,6,0,0,0)}|}+\frac{16(2,3,6,0,0,0)(M)}{|\mbox{Aut(2,3,6,0,0,0)}|},
\mylabel{eq:2.36}
\eeq
with $M\equiv 1\mymod 8$.
Analogously, \eqn{2.30} may be stated  as
\beq
\begin{split}
3(1,8,8,0,0,0)\left(\frac{M}{3^2}\right) & =(-n_1|3)\frac{16(1,6,6,0,0,0)(M)}{|\mbox{Aut(1,6,6,0,0,0)}|} \\
& +(-n_2|3)\frac{16(2,3,6,0,0,0)(M)}{|\mbox{Aut(2,3,6,0,0,0)}|},
\end{split}
\mylabel{eq:2.37}
\eeq
where $M\equiv 1\mymod 8$, $3|M$ and $n_1$, $n_2$ are any integers prime to $3$ that are represented 
by $x^2+6y^2+ 6z^2$, $2x^2+3y^2+6z^2$, respectively.

\bigskip
\section{Ramanujan's modular equations of degree $5$  and associated
identities for ternary quadratic forms with discriminant $400$} \label{sec:400}
\medskip

If all we ever wanted was to show that
\beq
\psi(q)(\phi(q)^2 -\phi(q^5)^2)\in P[q],
\mylabel{eq:3.1}
\eeq
we would be done in a second. Indeed, using an elementary trick that goes back to antiquity, we have
$5(x^2+y^2)=(x-2y)^2+(y+2x)^2$. Hence,
\beqs
\phi(q)^2-\phi(q^5)^2\in P[q],
\eeqs
and the truth of \eqn{3.1} becomes a bit boring. However, we want much more.
We ask for analogues of \eqn{2.18} and for associated theta-function identities.
Where do we begin? How about if we begin with binary forms of discriminant $-40$.
Again, there are just two genera of binary quadratic forms with this
discriminant: the class of $x^2+10y^2$ and the class of  $2x^2 +5y^2$.
We now add $10z^2$ to both forms to obtain ternaries $x^2+10y^2+10z^2$, $2x^2+5y^2+10z^2$ of discriminant $400$.
We observe that $2x^2+5y^2 +10z^2$ is the only form in its genus and that the
genus containing $x^2+ 10y^2+10z^2$ contains one more non-diagonal ternary form $4x^2+5y^2+6z^2+4zx$.

   It would be wrong to assume that we constructed all genera of ternary quadratic forms of discriminant $400$ this way.
In fact, we are being very selective by picking just two out of thirteen possible genera of discriminant $400$.
For the interested reader, we note that a table of genera of ternary quadratic forms, up to a large discriminant, is available on Neil Sloane's website at 

http://www.research.att.com/$\sim$njas/lattices/Brandt\_1.html\newline and was computed by Alexander Schiemann. In particular, this table includes relevant discriminants $144$, $400$, and $784$. The reader should be cautioned that the integer sextuple defining each form is preceded by an identification number and a colon, and that the identification number has no mathematical significance. The present authors use a combination of
Schiemann's software, scripts in a language called Magma, and C++ code written by the second author. For the reader with no experience of ternary forms we heartily recommend \cite{Dick}, especially the tables on pages 111-113.

  Again, it is easy to see that $(-n_1| 5)=1$ for any integer $n_1$ represented by the genus of $x^2+10y^2+10z^2$, 
$(n_1,5)=1$ and that $(-n_2|5 )= -1$ for any integer $n_2$ represented by $2x^2+5y^2+10z^2$, $(n_2,5) =1$.
Also 
$$|\mbox{Aut(1,10,10,0,0,0)}| =16,$$ 
$$|\mbox{Aut(4,5,6,0,4,0)}|=|\mbox{Aut(2,5,10,0,0,0)}|=8.$$
And so we anticipate two results similar to \eqn{2.36}, \eqn{2.37}. Namely,
\beq
\begin{split}
(1,8,8,0,0,0)(M)=(1,10,10,0,0,0)(M) & + 2(4,5,6,0,4,0)(M) \\
& + 2(2,5,10,0,0,0)(M),
\end{split}
\mylabel{eq:3.2}
\eeq
with $M\equiv 1\mymod 4$, and
\beq
\begin{split}
5(1,8,8,0,0,0)\left(\frac{M}{5^2}\right)=(1,10,10,0,0,0)(M) & + 2(4,5,6,0,4,0)(M) \\
& - 2(2,5,10,0,0,0)(M), 
\end{split}
\mylabel{eq:3.3}
\eeq
with $M\equiv 1\mymod 8$, $5|M$.
To prove \eqn{3.2} we rewrite it as
\beq
S_{8,1}(\phi(q)\phi(q^8)^2-\phi(q)\phi(q^{10})^{2}-2\phi(q^5)\chi(q)-2\phi(q^2)\phi(q^5)\phi(q^{10}))=0, 
\mylabel{eq:3.4}
\eeq
where
\beqs
\chi(q):= \sum_{x,z\in Z}q^{4x^2+4xz+6z^2}.
\eeqs
It is easy to see that
\beqs
\chi(q) = \sum_{x\equiv z\mymod 2} q^{x^2+5z^2}.
\eeqs
Hence,
\beq
\chi(q) = \phi(q^4)\phi(q^{20}) + 4q^6\psi(q^8)\psi(q^{40}). 
\mylabel{eq:3.5}
\eeq
Upon employing \eqn{2.4}, \eqn{3.5} we obtain 
\beqs
\begin{split}
S_{8,1}(\phi(q)\phi(q^8)^2) & = 2\psi(q)\phi(q)^2, \\
S_{8,1}(\phi(q)\phi(q^{10})^2) & = 2\psi(q)\phi(q^5)^2, \\
S_{8,1}(\phi(q^2)\phi(q^5)\phi(q^{10})) & = 8q^2\psi(q^2)\psi(q^5)\psi(q^{10}), \\
S_{8,1}(\phi(q^5)\chi(q)) & = S_{8,1}(\phi(q^5)\phi(q^4)\phi(q^{20})) \\
& = S_{8,1}(\phi(q^5)\phi(q^{16})\phi(q^{20})) + S_{8,1}(2q^4\phi(q^5)\psi(q^{32})\phi(q^{20})) \\
& = 4q^3\psi(q^5)\psi(q^{20})\phi(q^2)+4q\psi(q^4)\psi(q^5)\phi(q^{10}).
\end{split}
\eeqs
This means that \eqn{3.4} and, as a result, \eqn{3.2} is equivalent to the following theta-function identity
\beq
\begin{split}
\psi(q)(\phi(q)^2-\phi(q^5)^2) & =4q^3\psi(q^5)\psi(q^{20})\phi(q^2) \\
& +4q\psi(q^4)\psi(q^5)\phi(q^{10})+8q^2\psi(q^2)\psi(q^5)\psi(q^{10}).  
\end{split}
\mylabel{eq:3.6}
\eeq
It is easy to convert \eqn{3.6} into a modular equation of degree $5$,
\beq
\begin{split}
(m -1)\frac{\alpha^{\frac{1}{8}}}{\beta^{\frac{1}{8}}} & =(1+(1-\alpha)^{\frac{1}{2}})^{\frac{1}{2}}(1- (1-\beta)^{\frac{1}{2}})^{\frac{1}{2}} \\
& +(1-(1-\alpha)^{\frac{1}{2}})^{\frac{1}{2}}(1+(1-\beta)^{\frac{1}{2}})^{\frac{1}{2}}+2(\alpha\beta)^{\frac{1}{4}}.
\end{split}
\mylabel{eq:3.7}
\eeq
Here
\beqs
\begin{split}
m & =\frac{\phi(q)^2}{\phi(q^5)^2}, \\
\alpha & = 1-\frac{\phi(-q)^4}{\phi(q)^4}=16q\frac{\psi(q^2)^4}{\phi(q)^4}, \\
\beta & = 1-\frac{\phi(-q^5)^4}{\phi(q^5)^4}=16q^5\frac{\psi(q^{10})^4}{\phi(q^5)^4}.
\end{split}
\eeqs
Our proof of \eqn{3.7} hinges upon three powerful results established in \cite{Bern}, pp 285--286. 
\beq
2(1-(\alpha\beta)^{\frac{1}{2}}-((1-\alpha)(1-\beta))^{\frac{1}{2}})=(m-1)\left(-1+\frac{5}{m}\right), 
\mylabel{eq:3.8}
\eeq
\beq
\frac{\alpha^{\frac{1}{4}}}{\beta^{\frac{1}{4}}}= \frac{2m+r}{m(m-1)}, 
\mylabel{eq:3.9}
\eeq
\beq
4(\alpha^3\beta)^{\frac{1}{8}}=\frac{r}{m}+3-\frac{5}{m}, 
\mylabel{eq:3.10}
\eeq
where $r =(m(m^2-2m+5))^{\frac{1}{2}}$.
We begin by rewriting \eqn{3.7} as
\beq
\begin{split}
(m-1)\frac{\alpha^{\frac{1}{8}}}{\beta^{\frac{1}{8}}}-2(\alpha\beta)^{\frac{1}{4}} & =
(1+(1-\alpha)^{\frac{1}{2}})^{\frac{1}{2}}(1-(1-\beta)^{\frac{1}{2}})^{\frac{1}{2}} \\
& + (1-(1-\alpha)^{\frac{1}{2}})^{\frac{1}{2}}(1+(1-\beta)^{\frac{1}{2}})^{\frac{1}{2}}.
\end{split}
\mylabel{eq:3.11}
\eeq
Then we square both sides to obtain
\beqs
(m-1)^2\frac{\alpha^{\frac{1}{4}}}{\beta^{\frac{1}{4}}}-4(m-1)\alpha^{\frac{3}{8}}\beta^{\frac{1}{8}} =
2(1-(\alpha\beta)^{\frac{1}{2}}-((1-\alpha)(1-\beta))^{\frac{1}{2}}).
\eeqs
Next, we use \eqn{3.8}--\eqn{3.10} to arrive at the trivial statement
\beqs
\frac{(m-1)(2m+r)}{m}-\frac{(m-1)(r+3m-5)}{m}=(m-1)\left(\frac{5}{m}-1\right).
\eeqs
Hence, the proof of \eqn{3.7} is complete. Consequently, \eqn{3.2} is true, as desired.

To prove \eqn{3.3} we wish to consider  another modular equation of degree $5$
\beq
\begin{split}
\left(\frac{5}{m}-1\right)\frac{\beta^{\frac{1}{8}}}{\alpha^{\frac{1}{8}}} + 4(\alpha\beta)^{\frac{1}{4}} & = 
(1+(1-\alpha)^{\frac{1}{2}})^{\frac{1}{2}}(1-(1-\beta)^{\frac{1}{2}})^{\frac{1}{2}} \\
& +(1-(1-\alpha)^{\frac{1}{2}})^{\frac{1}{2}}(1+(1-\beta)^{\frac{1}{2}})^{\frac{1}{2}}.
\end{split}
\mylabel{eq:3.12}
\eeq
Comparing it with \eqn{3.11} we see that
\beqs
(m-1)\frac{\alpha^{\frac{1}{8}}}{\beta^{\frac{1}{8}}} =
\left(\frac{5}{m}-1\right)\frac{\beta^{\frac{1}{8}}}{\alpha^{\frac{1}{8}}}+ 4(\alpha\beta)^{\frac{1}{4}}.
\eeqs
Next, we multiply both sides by $\frac{\alpha^{\frac{1}{8}}}{\beta^{\frac{1}{8}}}$ to obtain 
\beqs
(m-1)\frac{\alpha^{\frac{1}{4}}}{\beta^{\frac{1}{4}}} = \frac{5}{m}-1+4(\alpha^3\beta)^{\frac{1}{8}}.
\eeqs
Employing \eqn{3.9}--\eqn{3.10}, we arrive at the trivial statement
\beqs
\frac{2m+r}{m} = \frac{5}{m}-1+\frac{r}{m}+3-\frac{5}{m}.
\eeqs
This completes the proof of \eqn{3.12}.
To proceed further we rewrite \eqn{3.12} in terms of theta-functions as
\beq
\begin{split}
5\psi(q^5)\phi(q^5)^2-\psi(q^5)\phi(q)^2 & = 4\psi(q)\psi(q^4)\phi(q^{10}) \\
& - 8q\psi(q)\psi(q^2)\psi(q^{10})+ 4q^2\psi(q)\psi(q^{20})\phi(q^2). 
\end{split}
\mylabel{eq:3.13}
\eeq
Using \eqn{2.4} and \eqn{3.5} and some elbow grease we can show that \eqn{3.13} is equivalent to
\beqs
S_{40,25}(\phi(q)\phi(q^{10})^2 + 2\phi(q^5)\chi(q)-2\phi(q^2)\phi(q^5)\phi(q^{10}))=10\psi(q^5)\phi(q^5)^2.
\eeqs
This implies that
\beqs
\begin{split}
& (1,10,10,0,0,0)(M)+2(4,5,6,0,4,0)(M)-2(2,5,10,0,0,0)(M) = \\
& 5\left|\left\{ n= 5x^2+ 5y^2 + 5\frac{(1+z)z}{2}:x,y,z\in\mathbf Z\right\}\right|, 
\end{split} 
\eeqs
with  $M=40n+25$.
The last equation can be easily recognized as \eqn{3.3}.
As before we can extend \eqn{3.2},\eqn{3.3} by using
$\frac{1}{3}(1,1,1,0,0,0)(M)$  in place of $(1,8,8,0,0,0)(M)$ as
\beq
\begin{split}
\frac{1}{3}(1,1,1,0,0,0)(M) & = (1,10,10,0,0,0)(M) \\
& + 2(4,5,6,0,4,0)(M)+2(2,5,10,0,0,0)(M),
\end{split} 
\mylabel{eq:3.14}
\eeq
with $M\equiv 1,2\mymod 4$, and
\beq
\begin{split}
\frac{5}{3}(1,1,1,0,0,0)\left(\frac{M}{5^2}\right) & =(1,10,10,0,0,0)(M) \\
& + 2(4,5,6,0,4,0)(M)-2(2,5,10,0,0,0)(M), 
\end{split} 
\mylabel{eq:3.15}
\eeq
with $M\equiv 1,2\mymod 4$, $5|M$.
While we have to suppress the details for the sake of brevity, we can not resist displaying four relevant theta-function identities.
\beqs
\begin{split}
\psi(q^2)(\phi(q)^2-\phi(q^5)^2)=2q\psi(q^5)^2\phi(q) & + 2q\psi(q^{10})\phi(q^2)\phi(q^{10}) \\
& + 8q^4\psi(q^4)\psi(q^{10}\psi(q^{20}).
\end{split}
\eeqs
This one proves the $M\equiv 1\mymod 4$ case in \eqn{3.14}.
Analogously,
\beqs
\begin{split}
\psi(q^2)^2\phi(q)& = q^2\psi(q^{10})^2\phi(q)+2q\psi(q^2)\psi(q^5)^2 \\
& + q^2\psi(q^{20})\phi(q^2)\phi(q^5)+\psi(q^4)\phi(q^5)\phi(q^{10})
\end{split}
\eeqs
proves the $M\equiv 2\mymod 4$ case in \eqn{3.14}.
Finally,
\beqs
\begin{split}
5q\psi(q^{10})\phi(q^5)^2 & = q\psi(q^{10})\phi(q)^2 \\
& + 2\psi(q)^2\phi(q^5)-8q^3\psi(q^2)\psi(q^4)\psi(q^{20})-2\psi(q^2)\phi(q^2)\phi(q^{10})
\end{split}
\eeqs
and
\beqs
\begin{split}
5q^2\psi(q^5)^2\psi(q^{10}) & = \psi(q^2)^2\phi(q^5) \\
& + 2q\psi(q)^2\psi(q^{10})-\psi(q^4)\phi(q)\phi(q^{10})-q^2\psi(q^{20})\phi(q)\phi(q^2)
\end{split}
\eeqs
are required to prove \eqn{3.15}.

\bigskip
\section{Ternary forms with discriminant $784$.} \label{sec:784}
\medskip
Here we will prove the Theorem 1 and Theorem 2, stated in the Introduction.
It seems that the identities for theta functions in \eqn{1.14}, \eqn{1.15}
correspond to modular equations of mixed degree $21$.  While Ramanujan had some results for modular
equations of this degree \cite{Bern}, we could not find enough relations to handle our formulas \eqn{1.14}, \eqn{1.15}.
And so it is with some reluctance that we resort to routine modular function techniques. The necessary
background theory on modular functions and forms may be found in Rankin's book \cite{Ran}. Of central
importance to us is the valence formula (p.98,  \cite{Ran}).

We begin by dividing  both sides of \eqn{1.14}  by  $\psi(q)\phi(q)^2$.
Making use of \eqn{1.7}, \eqn{1.8}, \eqn{1.9}, \eqn{1.10}, \eqn{1.12} we end up with a simple identity for four $\eta$-quotients
\beq
g_{1}(z)+4g_{2}(z)+8g_{3}(z)+8g_{4}(z)=1,
\mylabel{eq:4.1}
\eeq
where
\beqs
g_{1}(z): = \frac{\eta(14z)^{10}\eta(4z)^4\eta(z)^4}{\eta(28z)^4\eta(7z)^4\eta(2z)^{10}},
\eeqs
\beqs
g_{2}(z): = \frac{\eta(14z)^7\eta(4z)^6\eta(z)^5}{\eta(28z)^2\eta(7z)^3\eta(2z)^{13}},
\eeqs
\beqs
g_{3}(z): = \frac{\eta(42z)^5\eta(28z)^2\eta(6z)^2\eta(4z)^4\eta(z)^5}{\eta(84z)^2\eta(21z)^2\eta(14z)\eta(3z)\eta(2z)^{12}},
\eeqs
\beqs
g_{4}(z): = \frac{\eta(84z)\eta(28z)\eta(21z)\eta(14z)\eta(4z)^4\eta(3z)^2\eta(z)^4}
{\eta(42z)\eta(7z)\eta(6z)\eta(2z)^{11}}.
\eeqs
Clearly, all our $\eta$-quotients in \eqn{4.1} are of the form
\beqs
f(z) = \prod_{\delta|n }\eta(\delta z)^{r_{\delta}},
\eeqs
where $n$ is some positive integer ($84$  in our case) and all $ \delta\ge 1$, $r_{\delta}$ are integers.
The following result was proved by Morris Newman \cite{New}.
\begin{theorem}\label{N}
The $\eta$-quotient $f(z)$ is a modular function on
\beqs
\Gamma_{0}(n): = \bigg\{\begin{bmatrix} a & b  \\ f & c 
\end{bmatrix}
\in SL_2(Z):c\equiv 0\mymod n\bigg\},
\eeqs
if the following four conditions are met
\begin{align*}
& \sum_{\delta|n }r_{\delta}=0, \\
& \sum_{\delta|n}\delta r_{\delta}\equiv 0\mymod 24, \\
& \sum_{\delta|n}\frac{nr_{\delta}}{\delta}\equiv 0\mymod 24, \\
& \prod_{\delta|n }\delta^{r_{\delta}} \quad\mbox{is a rational square.}
\end{align*}

\end{theorem}
It is now straightforward to verify that $g_{1}(z),g_{2}(z),g_{3}(z),g_{4}(z)$
are modular functions on $\Gamma_{0}(84)$.
Consequently,
\beqs
h(z): = g_{1}(z)+4g_{2}(z)+8g_{3}(z)+8g_{4}(z)-1
\eeqs
is also a modular function on $\Gamma_{0}(84)$.
To proceed we will need the following observation from \cite{Bia}
\begin{theorem}\label{B}
If $n$ is square free integer, then a complete set of inequivalent
cusps for $\Gamma_{0}(4n)$ is $\{\frac{1}{s}:s |4n \}$.
\end{theorem}
And so $k\cup\{\frac{1}{84}\}$ is a complete set of $12$ inequivalent cusps of
$\Gamma_0(84)$ where  $k:=\{1,\frac{1}{2},\frac{1}{6},\frac{1}{4},\frac{1}{12},\frac{1}{7},\frac{1}{42},\frac{1}{21},\frac{1}{3},\frac{1}{14},\frac{1}{28}\}$.
From the definition of $\eta(z)$, it follows that $\eta$-quotients have no zeros or poles in the upper-half plane 
(i.e. $Im(z)> 0$). Ligozat \cite{Lig} calculated the order of $\eta$-quotient $f(z)$ at the cusps of $\Gamma_0(n)$.
\begin{theorem}\label{L} 
If an $\eta$-quotient $f$ is a modular function on $\Gamma_{0}(n)$, then at the cusp $\frac{b}{c}$ with $\gcd(b,c)=1$
\beqs
ORD\left(f,\frac{b}{c}\right)=\frac{n}{24\gcd(n,c^2)}\sum_{\delta|n}\frac{r_{\delta}\gcd(c,\delta)^2}{\delta}.
\eeqs
\end{theorem}
This way we obtain the following

\medskip
\centerline{TABLE 1}
\[
\arraycolsep = 7mm
\begin{array}{crrrrr}
\mbox{CUSP}  &     O_{1}(s)  &  O_{2}(s) &   O_{3}(s)  & O_{4}(s)  &  O_{h}(s) \\
  1   &    0   &   0  &    0  &    0  &    0 \\
 1/2  &   -9  &  -12  &  -12  &  -12  &  -12 \\
 1/6  &   -3   &  -4  &   -1  &   -4  &   -4 \\
 1/4  &    0   &   3  &   -1  &   -1  &   -1 \\
1/12  &    0   &   1   &   2   &   0  &    0 \\
 1/7  &    0   &   0   &   0   &   0   &   0 \\
1/42  &    3   &   2   &   5   &   0   &   0 \\
1/21  &    0   &   0   &   0  &    3   &   0 \\
 1/3  &    0   &   0   &   0   &   5   &   0 \\
1/14  &    9   &   6   &   0   &   0   &   0 \\
1/28  &    0   &   3   &   5   &   5   &   0 \\
\end{array}
\]
where
$O_{h}(s)$ is lower bound for $ORD(h,s)$ and $O_{i}(s):= ORD(g_{i},s), i=1,2,3,4$.
To prove \eqn{1.14} and \eqn{4.1} we must show that $h(z)=0$.
The valence formula implies that (unless $h$ is a constant)
\beqs
\sum_{s\in k}ORD(h,s)+ORD(h,1/84)\leq 0.
\eeqs
Using data collected in Table 1 and keeping in mind that cusp $\frac{1}{84}$ is equivalent to $i\infty$, we infer that  (unless $h$ is constant)
\beqs
-17+ORD(h,i\infty)\leq 0.
\eeqs
But direct inspection shows that $ORD(h,i\infty) > 17$.
That is, if one expands $h$ in powers of $q$, then one finds that the first $18$ coefficients in this expansion are zero.
Hence, one arrives at a contradiction. This contradiction implies that $h=0$, as desired.
This completes our proof of \eqn{1.14}.
Obviously, we could've proved \eqn{1.15} in a similar fashion. Instead, we choose a more painful way,
because there is nothing like pain for achieving excellence.
In any event, our  approach will shed some extra light on the relation between \eqn{1.14} and \eqn{1.15}.
It is not hard  to verify (in term by term fashion) that
\beq
\begin{split}
& 8q^2\psi(q^3)\psi(q^{14})\phi(q^{21})+  8q^4\psi(q^{14})f(q,q^2)f(q^7,q^{35}) = \\
& C(q)(8q^2\psi(q^2)\psi(-q^{21})\phi(-q^3)+8\psi(q^2)f(-q,-q^5)f(-q^7,q^{14})), 
\end{split}
\mylabel{eq:4.2}
\eeq
where
\beqs
C(q)=q^2\frac{E(q^{28})E(q^{14})E(q^2)^2}{E(q^7)E(q^4)^2E(q)}.
\eeqs
All one needs is a simple formula
\beqs
E(-q)=\frac{E(q^2)^3}{ E(q^4)E(q)}.
\eeqs
We will prove shortly that
\beq
\begin{split}
& \psi(q)(\phi(q)^2-\phi(q^7)^2)-4q\psi(q^2)\psi(q^7)\phi(q^7)= \\
& C(q)(\psi(-q^7)(7\phi(-q^7)^2+\phi(-q)^2)-4q\psi(-q)\psi(q^{14})\phi(-q)).
\end{split}
\mylabel{eq:4.3}
\eeq
Next, we rewrite \eqn{1.14} as
\beq
\begin{split}
& \psi(q)(\phi(q)^2-\phi(q^7)^2)-4q\psi(q^2)\psi(q^7)\phi(q^7)= \\
& 8q^2\psi(q^{3})\psi(q^{14})\phi(q^{21})+8q^4\psi(q^{14})f(q,q^2)f(q^7,q^{35}),
\end{split}
\mylabel{eq:4.4}
\eeq
and use \eqn{4.2} on the right and \eqn{4.3} on the left to get
\beq
\begin{split}
& C(q)(\psi(-q^7)(7\phi(-q^7)^2+\phi(-q)^2)-4q\psi(-q)\psi(q^{14})\phi(-q))= \\
& C(q)(8q^2\psi(q^2)\psi(-q^{21})\phi(-q^3)+8\psi(q^2)f(-q,-q^5)f(-q^7,q^{14})).
\end{split}
\mylabel{eq:4.5}
\eeq
Dividing both sides by $C(q)$ and replacing $q$ by $-q$, we get \eqn{1.15}.

But what about \eqn{4.3}? We start by rewriting it as a modular equation of degree $7$.
Namely,
\beq
\begin{split}
m-1-2t=\frac{7}{mt^2}\frac{\beta^{\frac{7}{12}}(1-\beta)^{\frac{7}{12}}}{\alpha^{\frac{1}{12}}(1-\alpha)^{\frac{1}{12}}}
& + \frac{1}{t^2}\frac{\beta^{\frac{7}{12}}(1-\beta)^{\frac{7}{12}}}{\alpha^{\frac{1}{12}}(1-\alpha)^{\frac{1}{12}}}            \frac{(1-\alpha)^{\frac{1}{2}}}{(1-\beta)^{\frac{1}{2}}} \\
& - 2\frac{\beta}{t^4} \frac{\alpha^{\frac{7}{24}}(1-\alpha)^{\frac{7}{24}}}{\beta^{\frac{1}{24}}(1-\beta)^{\frac{1}{24}}},
\end{split}
\mylabel{eq:4.6}
\eeq
Here $m$ is a multiplier of degree seven, $\beta$ has degree seven over $\alpha$ and $t$ is defined by
\beqs
t=(\alpha\beta)^{1/8}.
\eeqs
 The proof of the following modular equations of degree seven can be found in (p.314, \cite{Bern})
\beq
(1-\alpha)^{\frac{1}{8}}(1-\beta)^{\frac{1}{8}} = 1-t, 
\mylabel{eq:4.7}
\eeq
\beq
\frac{7(1-2t)}{m}+1=4\frac{\alpha^{\frac{7}{24}}(1-\alpha)^{\frac{7}{24}}}{\beta^{\frac{1}{24}}(1-\beta)^{\frac{1}{24}}},
\mylabel{eq:4.8}
\eeq
\beq
\frac{(m(2t-1)+1)^{2}}{16}=\frac{\beta^{\frac{7}{12}}(1-\beta)^{\frac{7}{12}}}{\alpha^{\frac{1}{12}}
(1-\alpha)^{\frac{1}{12}}}.
\mylabel{eq:4.9}
\eeq
Observe that \eqn{4.7} implies that
\beqs
\frac{(1-\alpha)^{\frac{1}{2}}}{(1-\beta)^{\frac{1}{2}}}= \frac{(1-t)^4}{1-\beta}.
\eeqs
And so \eqn{4.6} becomes
\beq
m-1-2t+\beta\frac{\frac{7(1-2t)}{m}+1 }{ 2t^4}-\left(\frac{7}{m}+ \frac{(1-t)^4}{1-\beta}\right)\frac{(m(2t-1)+1)^2}{16t^2}=0.
\mylabel{eq:4.10}
\eeq
Next, we recall the equation (19.19) in \cite{Bern}
\beqs
m=\frac{t-\beta}{t(1-t)(1-t+t^2)}.
\eeqs
We use it to eliminate $m$ from \eqn{4.10} and obtain after some algebra that
\beq
(\beta^2-\beta((1+t^8)-(1-t)^8)+t^8)\frac{P(t,\beta)}{Q(t,\beta)}=0, 
\mylabel{eq:4.11}
\eeq
where $P(t,\beta)$ and $Q(t,\beta)$ are some polynomials in $t$ and $\beta$.
But $\beta$ and $\alpha$ are roots of the quadratic equation
\beqs
x^2-x((1+t^8)-(1-t)^8)+t^8 =0,
\eeqs
as observed on page 316 in \cite{Bern}. Hence, the proof of \eqn{4.3} is complete.
Consequently, \eqn{1.15} is true.

We now ready to move on to prove our Theorem 2.
Clearly, \eqn{1.16} is equivalent to
\beq
S_{8,1}\left(\phi(q)\phi(q^8)^2-\phi(q)\phi(q^{14})^2-2\phi(q^2)\phi(q^7)\phi(q^{14})-4\phi(q^{14})u(q)\right)=0, 
\mylabel{eq:4.12}
\eeq
where
\beqs
u(q):=\sum_{x,z\in Z}q^{3x^2+2xy+ 5y^2}.
\eeqs
It is not hard to check that
\beqs
u(q)=\sum_{x\equiv y\mymod 3}q^{\frac{x^2+ 14y^2}{3}}.
\eeqs
Hence,
\beq
u(q)=\phi(q^3)\phi(q^{42})+2q^5f(q,q^5)f(q^{14},q^{70}).  
\mylabel{eq:4.13}
\eeq
We shall also require
\beq
S_{8,1}(q^3f(q,q^5)f(q^{14},q^{70}))=q^2f(q,q^2)f(q^7,q^{35}). 
\mylabel{eq:4.14}
\eeq
Using \eqn{2.4} together with \eqn{4.13} and \eqn{4.14}, we deduce that
\beq
\begin{split}
S_{8,1}(\phi(q^{14})u(q )) & = 2q\psi(q^{14})S_{8,1}( q^6(\phi(q^3)\phi(q^{42})+2q^5 f(q,q^5)f(q^{14},q^{70}))) \\
& = 4q^2\psi(q^{14})\psi(q^3)\phi(q^{21})+4q^4\psi(q^{14})f(q,q^2)f(q^7,q^{35}).
\end{split}
\mylabel{eq:4.15}
\eeq
Next, with the aid of \eqn{2.4} we verify that
\beq
S_{8,1}(\phi(q)\phi(q^{14})^2)=2\psi(q)\phi(q^7)^2,  
\mylabel{eq:4.16}
\eeq
\beq
S_{8,1}(\phi(q^2)\phi(q^7)\phi(q^{14}))=4q\psi(q^2)\psi(q^7)\phi(q^7).
\mylabel{eq:4.17}
\eeq
Combining \eqn{2.14}, \eqn{4.12} and \eqn{4.15}--\eqn{4.17} we end up with \eqn{1.14}.Hence, the proof of \eqn{1.16} is complete. 

Our proof of \eqn{1.17} is analogous.
We verify that
\beqs
S_{56,49}(\phi(q)\phi(q^{14})^2)=2\psi(q^7)\phi(q)^2,
\eeqs
\beqs
S_{56,49}(\phi(q^2)\phi(q^7)\phi(q^{14}))=4q\psi(q)\psi(q^{14})\phi(q),
\eeqs
\beqs
S_{56,49}(\phi(q^{14})u(q ))=4q^2\psi(q^2)\psi(q^{21})\phi(q^3)+4\psi(q^2)f(q,q^5)f(q^7,q^{14}).
\eeqs
These results enable us to convert \eqn{1.15} into
\beqs
S_{56,49}\left(-\phi(q)\phi(q^{14})^2-2\phi(q^2)\phi(q^7)\phi(q^{14})+4\phi(q^{14})u(q)\right)=14\psi(q^7)\phi(q^7)^2.
\eeqs
Hence, we have 
\beqs
\begin{split}
& - (1,14,14,0,0,0)(M)-2(2,7,14,0,0,0)(M)+4(3,5,14,0,0,2)(M) \\
& = 7\left|\left\{n = 7\frac{x(x+1)}{2}+7y^2+7z^2:x,y,z\in\mathbf Z\right\}\right| \\
& = 7\left|\left\{1+8\frac{n}{7}= x^2+8y^2+8z^2:x,y,z\in\mathbf Z\right\}\right|,\quad M=56n+49.
\end{split}
\eeqs
The truth of \eqn{1.17} is now transparent.

Again, we can extend \eqn{1.14},\eqn{1.15} by using $\frac{1}{3}(1,1,1,0,0,0)(M)$ instead of $(1,8,8,0,0,0)(M)$.
This way we have  
\beq
\begin{split}
\frac{1}{3}(1,1,1,0,0,0)(M)=(1,14,14,0,0,0)(M) & + 2(2,7,14,0,0,0)(M) \\
& + 4(3,5,14,0,0,2)(M),
\end{split}
\mylabel{eq:4.18}
\eeq
with $M\equiv 1,2\mymod 4$, and
\beq
\begin{split}
\frac{7}{3}(1,1,1,0,0,0)\left(\frac{M}{7^2}\right)=-(1,14,14,0,0,0)(M) & - 2(2,7,14,0,0,0)(M) \\
& + 4(3,5,14,0,0,2)(M),
\end{split}  
\mylabel{eq:4.19}
\eeq
with $M\equiv 1,2\mymod 4$, $7|M$.
We limit ourselves to a few remarks. The theta-function identities
\beqs
\begin{split}
\psi(q^2)(\phi(q)^2-\phi(q^7)^2)& = 4q^2\psi(q^4)\psi(q^{14})\phi(q^{14})+4q^5\phi(q^2)\psi(q^{14})\psi(q^{28}) \\
& + 8q^4\psi(q^6)\psi(q^7)\psi(q^{21})+4q\psi(q^7) f(q^2,q^4)f(q^7,q^{14}),
\end{split}  
\eeqs
and
\beqs
\begin{split}
\phi(q)(\psi(q^2)^2-q^3\psi(q^{14})^2)& = q^3\psi(q^{28})\phi(q^2)\phi(q^7)+\psi(q^4)\phi(q^7)\phi(q^{14}) \\
& + 2q^3\psi(q^7)\psi(q^{21})\phi(q^3)+2q\psi(q^7)f(q,q^5)f(q^7,q^{14})
\end{split}  
\eeqs
imply the $M\equiv 1\mymod 4$ and $M\equiv2\mymod 4$ cases in \eqn{4.18}, respectively.
Analogously,
\beqs
\begin{split}
7q\psi(q^7)^2\phi(q^7) & = 8q^5\psi(q)\psi(q^3)\psi(q^{42})+4\psi(q)f(q,q^2)f(q^{14},q^{28}) \\
& - q\phi(q)^2\psi(q^{14})-4\psi(q^2)\psi(q^4)\phi(q^{14})-4q^3\psi(q^2)\psi(q^{28})\phi(q^2)
\end{split}  
\eeqs
and
\beqs
\begin{split}
7q^3\psi(q^7)^2\psi(q^{14})& = 2\psi(q)\psi(q^3)\phi(q^{21})+2q^2\psi(q)f(q^7,q^{35})f(q,q^2) \\
& - \psi(q^2)^2\phi(q^7)-\psi(q^4)\phi(q)\phi(q^{14})-q^3\psi(q^{28})\phi(q)\phi(q^2) 
\end{split}  
\eeqs
can be used to prove both cases in \eqn{4.19}.

We conclude this section by showing how to relate the ternaries on the right of \eqn{1.16} and \eqn{1.17} to binaries 
with discriminant $-56$. Again, there are just two genera of binary quadratic forms with this discriminant.
The first one contains the class of $x^2+14y^2$ and the class of $2x^2+7y^2$. The second one contains the class of 
$3x^2+2xy+5y^2$ and the class of $3x^2-2xy+5y^2$. We now add $14z^2$ to both forms in the first genus of binary
quadratic forms to obtain our first genus of ternary quadratic forms of discriminant $784$
\beqs
\{x^2+14y^2+14z^2, \quad 2x^2+7y^2+14z^2\}.
\eeqs
Next, we add $14z^2$ to the first forms in the second genus of binary quadratic forms to obtain our second genus of ternary quadratic forms of discriminant $784$

\beqs
\{3x^2+2xy+5y^2+14z^2\}.
\eeqs

 Both these genera are complete as described, no other forms need be added. As a result,  $3 x^2 + 2 x y + 5 y^2 + 14 z^2$
  is obviously regular. Also, it is easy to verify that
  
$$\frac{16}{|\mbox{Aut}(1,14,14,0,0,0)|}=1,$$ 
$$\frac{16}{|\mbox{Aut}(2,7,14,0,0,0)|}=2,$$
$$\frac{16}{|\mbox{Aut}(3,5,14,0,0,2)|}=4.$$
And so \eqn{1.16} can be stated as
\beqs
\begin{split}
(1,8,8,0,0,0)(M)=\frac{16(1,14,14,0,0,0)(M)}{|\mbox{Aut}(1,14,14,0,0,0)|}
& + \frac{16(2,7,14,0,0,0)(M)}{|\mbox{Aut}(2,7,14,0,0,0)|} \\
& + \frac{16(3,5,14,0,0,2)(M)}{|\mbox{Aut}(3,5,14,0,0,2)|},
\end{split}
\eeqs
with $M\equiv 1\mymod 8$.
This shape of identity clearly demonstrates that on the right we have Siegel's type weighted average over two genera 
of ternary quadratic form. Moreover, it is easy to see that
\beqs
(-{n_{1}}|7)=-1,
\eeqs
for any integer $n_{1}$ represented by the genus of $x^2+14y^2+14z^2$, $(n_{1},7)=1$ and that
\beqs
(-{n_{2}}|7)=1,
\eeqs
for any integer $n_{2}$ represented by $3x^2+5y^2+14z^2+2xy$, $(n_{2},7)=1$. 
This allows us to rewrite \eqn{1.17} as
\beqs
\begin{split}
7(1,8,8,0,0,0)\left(\frac{M}{7^2}\right) & =(-n_{1}|7)\left(\frac{16(1,14,14,0,0,0)(M)}{|\mbox{Aut}(1,14,14,0,0,0)|}
+\frac{16(2,7,14,0,0,0)(M)}{|\mbox{Aut}(2,7,14,0,0,0)|}\right)\\
& + (-n_{2}|7)\frac{16(3,5,14,0,0,2)(M)}{|\mbox{Aut}(3,5,14,0,0,2)|},
\end{split}
\eeqs
with $M\equiv1\mymod 8$,$7|M$.
  
\bigskip
\section{Ternary forms with discriminant $3600$.} \label{sec:3600}
\medskip

Up to now all our theorems involved certain ternary forms with disciminant $16p^2$ for prime $p= 3,5,7$.  
In this section we consider ternaries with discriminant $16S^2$ for composite $S=15$. This case has all of the ingredients of the general case to be discussed later. Again, we start with the binaries with discriminant $-120$.
There are 4 genera with this discriminant and each has a single class per
genus $\{x^2+ 30y^2 \}$, $\{3x^2+10y^2\}$, $\{5x^2+6y^2\}$, $\{2x^2+ 15y^2\}$. Following a well-trodden path we add $30z^2$ to each of these forms to get four ternary forms with discriminant $3600$. Next, we extend each ternary form to a genus of ternary quadratic forms. This way we obtain four genera of ternary quadratic forms with discriminant $3600$:
$$
\{x^2+30y^2\}\rightarrow TG_1:=\{x^2+30y^2+30z^2,6x^2 +10y^2 +15z^2\},
$$
$$
\{3x^2+10y^2\}\rightarrow TG_2:=\{3x^2+10y^2+30z^2\},
$$
$$
\{5x^2+6y^2\}\rightarrow TG_3:=\{5x^2+6y^2+30z^2,9x^2 +11y^2+11z^2+2yz+6zx+6xy\},
$$
$$
\{2x^2+15y^2\}\rightarrow TG_4:=\{2x^2+15y^2 +30z^2,5x^2+12y^2+18z^2+12yz\}.
$$
We check that
$$
\frac{16}{|\mbox{Aut}(1,30,30,0,0,0)|}=1,
$$ 
$$
\frac{16}{|\mbox{Aut}(6,10,15,0,0,0)|}=\frac{16}{|\mbox{Aut}(3,10,30,0,0,0)|}=2
$$
$$
\frac{16}{|\mbox{Aut}(5,6,30,0,0,0)|}=\frac{16}{|\mbox{Aut}(2,15,30,0,0,0)|}=\frac{16}{|\mbox{Aut}(5,12,18,12,0,0)|}=2,
$$
$$
\frac{16}{|\mbox{Aut}(9,11,11,2,6,6)|}=4.
$$
We now take Siegel's weighted average over the four genera above. This way we are led to
\beq
\begin{split}
(1,8,8,0,0,0 )(M) & = (1,30,30,0,0,0)(M)+2(6,10,15,0,0,0)(M) \\
& + 2(3,10,30,0,0,0)(M)+2(5,6,30,0,0,0)(M) \\
& + 4(9,11,11,2,6,6)(M)+2(2,15,30,0,0,0)(M) \\
& + 2(5,12,18,12,0,0,0)(M), 
\end{split}
\mylabel{eq:5.1}
\eeq
with $M\equiv1\mymod 8$.
This can be stated compactly as
\beq
(1,8,8,0,0,0)(M)=\sum_{i=1}^{4}W_{i}(M), 
\mylabel{eq:5.2}
\eeq
with $M\equiv1\mymod 8$. 
Here,    
\beqs
W_{i}(M)=16\sum_{f\in TG_{i}} \frac{R_{f}(M)}{|\mbox{Aut}(f)|},\quad i=1,2,3,4,
\eeqs
and $R_{f}(M)$ denotes the number of representations of $M$ by $f$.
The associated theta-function identity is as follows
\beq
\begin{split}
\psi(q)\phi(q)^2 & = \psi(q)\phi(q^{15})^2+4q^3\psi(q^{10})\psi(q^{15})\phi(q^3)+4q^4\psi(q^3)\psi(q^{30})\phi(q^5) \\
& + 8q^5\psi(q^5)\psi(q^6)\psi(q^{30})+4q\psi(q^9)\phi(q^{45})^2 \\
& + 8q^{11}\psi(q^9)f(q^{15},q^{75})^2+8q^5\phi(q^{45})f(q^3,q^6)f(q^{15},q^{75}) \\
& + 4q^{10}f(q^3,q^6)f(q^{15},q^{75})^2+4q^2\psi(q^2)\psi(q^{15})\phi(q^{15}) \\
& + 4q^8\psi(q^5)\psi(q^{60})\phi(q^6)+4q^2\psi(q^5)\psi(q^{12})\phi(q^{30}).
\end{split}
\mylabel{eq:5.3}
\eeq
Using \eqn{5.3} we easily  see that
\beqs
\psi(q)(\phi(q)^2-\phi(q^{15})^2)\in P[q].
\eeqs
To prove \eqn{5.3} we divide both sides by $\psi(q)\phi(q)^2$ and use \eqn{1.7}--\eqn{1.10},
\eqn{1.12} to end up with an identity for eleven $\eta$-quotients

\beq
\begin{split}
1 & = \frac{\eta(30z)^{10}\eta(4z)^4\eta(z)^4}{\eta(60z)^4\eta(15z)^4\eta(2z)^{10}}
+4\frac{\eta(30z)^2\eta(20z)^2\eta(6z)^5\eta(4z)^4\eta(z)^5}{\eta(15z)\eta(12z)^2\eta(10z)\eta(3z)^2\eta(2z)^{12}} \\
& + 4\frac{\eta(60z)^2\eta(10z)^5\eta(6z)^2\eta(4z)^4\eta(z)^5}{\eta(30z)\eta(20z)^2\eta(5z)^2\eta(3z)\eta(2z)^{12}}
+8\frac{\eta(60z)^2\eta(12z)^2\eta(10z)^2\eta(4z)^4\eta(z)^5}{\eta(30z)\eta(6z)\eta(5z)\eta(2z)^{12}}\\
& + 4\frac{\eta(90z)^{10}\eta(18z)^2\eta(4z)^4\eta(z)^5}{\eta(180z)^4\eta(45z)^4\eta(9z)\eta(2z)^{12}}
+8\frac{\eta(180z)^2\eta(45z)^2\eta(30z)^4\eta(18z)^2\eta(4z)^4\eta(z)^5}
{\eta(90z)^2\eta(60z)^2\eta(15z)^2\eta(9z)\eta(2z)^{12}}\\
& + 8\frac{\eta(90z)^4\eta(30z)^2\eta(9z)^2\eta(6z)\eta(4z)^4\eta(z)^5}
{\eta(180z)\eta(60z)\eta(45z)\eta(18z)\eta(15z)\eta(3z)\eta(2z)^{12}} \\
& + 4\frac{\eta(180z)^2\eta(45z)^2\eta(30z)^4\eta(9z)^2\eta(6z)\eta(4z)^4\eta(z)^5}
{\eta(90z)^2\eta(60z)^2\eta(18z)\eta(15z)^2\eta(3z)\eta(2z)^{12}}\\
& + 4\frac{\eta(30z)^7\eta(4z)^6\eta(z)^5}{\eta(60z)^2\eta(15z)^3\eta(2z)^{13}}
+4\frac{\eta(120z)^2\eta(12z)^5\eta(10z)^2\eta(4z)^4\eta(z)^5}{\eta(60z)\eta(24z)^2\eta(6z)^2\eta(5z)\eta(2z)^{12}} \\
& + 4\frac{\eta(60z)^5\eta(24z)^2\eta(10z)^2\eta(4z)^4\eta(z)^5}{\eta(120z)^2\eta(30z)^2\eta(12z)\eta(5z)\eta(2z)^{12}}. 
\end{split}
\mylabel{eq:5.4}
\eeq
To verify the last identity, we use the Newman theorem stated in the last section to show that all eleven quotients 
on the right of \eqn{5.4} are modular functions on $\Gamma_{0}(360)$.
Next, we  confirm that the right hand side expanded in powers of $q$ is just $1$ for sufficient number of terms. 
Just how many terms to calculate is determined by the Ligozat theorem. In this regard we note that $K\cup\frac{1}{360}$
is a complete set of $32$ inequivalent cusps of $\Gamma_{0}(360)$.
Here,
\beqs
\begin{split}
K:= & \left\{1,\frac{1}{2},\frac{1}{3},\frac{2}{3},\frac{1}{4},\frac{1}{5},\frac{1}{6},\frac{5}{6},\frac{1}{8},\frac{1}{9}, \frac{1}{10},\frac{1}{12},\frac{5}{12},\frac{1}{15},\frac{2}{15},\frac{1}{18},\frac{1}{20},\frac{1}{24},\right. \\
& \left. \frac{5}{24},\frac{1}{30},\frac{11}{30},\frac{1}{36},\frac{1}{40},\frac{1}{45}, \frac{1}{60},\frac{11}{60}, \frac{1}{72},\frac{1}{90},\frac{1}{120},\frac{5}{120},\frac{1}{180}\right\}.
\end{split}
\eeqs
Let $H$ denote $-1+$ RHS \eqn{5.4}. Obviously, $H$ is a modular function on $\Gamma_{0}(360)$.
From Ligozat's theorem we see that

\medskip
\centerline{TABLE 2}
\[
 \arraycolsep = 2.77mm
\begin{array}{cccccccccccc}
\mbox{CUSP}  & 1 &   \frac{1}{2} &  \frac{1}{3} &  \frac{2}{3}  & 
\frac{1}{4}  &   \frac{1}{5}  &  \frac{1}{6}  &    \frac{5}{6} &   
\frac{1}{8}  &  \frac{1}{9}  \\
O_H(s)  & 0 &  -54 &    0 &     0 &     -5 &      0 &    -6 &  -6    &     -5
&     0  \\
   \\
\mbox{CUSP}  & \frac{1}{10}  &  \frac{1}{12}
 &   \frac{5}{12}  &  \frac{1}{15}  &  \frac{2}{15}  &  \frac{1}{18} &
\frac{1}{20}  &   \frac{1}{24} &   \frac{5}{24}  &   \frac{1}{30}  \\
O_H(s) &     -6 &      0 &  0 &       0 &          0 &      -6 &     -1 &    
  0 &       0 &        0 \\
        \\
 \mbox{CUSP}   &   \frac{11}{30}  &   \frac{1}{36} &   \frac{1}{40} &  
\frac{1}{45} &   \frac{1}{60}  &
  \frac{11}{60}  &   \frac{1}{72}  &   \frac{1}{90}   &   \frac{1}{120}  &
  \frac{5}{120}   &   \frac{1}{180}  \\
 O_H(s)  &     0  &     0 &  -1 &    0 &     0 &    0 &    0 &     0 &     0
&      0 &       0
\end{array}
\]
where $O_{H}(s)$ is lower bound for $ORD( H,s)$. The valence formula says that (unless $H$ is a constant)
\beqs
\sum_{s \in K }ORD(H,s)+ORD (H,\frac{1}{360}) \leq 0 .
\eeqs
Using data collected in the Table 2 and keeping in mind that cusp $\frac{1}{360}$
is equivalent to $i\infty$, we deduce that (unless $H$ is constant)
\beq
-90+ORD(H,i\infty)\leq 0. 
\mylabel{eq:5.5}
\eeq
We use Maple to calculate $91$ coefficients of the Fourier expansion of $H$.
This way we see that $ORD(H,i\infty)\geq 90$, and thus \eqn{5.5} is contradicted.
Hence, $H=0$. The proof of \eqn{5.3} and \eqn{5.4} is now complete.
In exactly the same mechanical manner we can prove three companion identities.
\beq
\begin{split}
3\psi(q^3)\phi(q^3)^2 & = 4q\psi(q^5)\psi(q^6)\phi(q^5)+8q^3\psi(q^2)\psi(q^{10})\psi(q^{15}) \\
& + 4\psi(q^3)\phi(q^{15})^2+4q^3f(q,q^2)f(q^5,q^{25})^2+4q^2\psi(q^4)\psi(q^{15})\phi(q^{10}) \\ 
& + 4q^4\psi(q^{15})\psi(q^{20})\phi(q^2)-\psi(q^3)\phi(q^5)^2-4q^4\psi(q^5)\psi(q^{30})\phi(q) \\
& - 4q\psi(q)\psi(q^{10})\phi(q^{15}),
\end{split}
\mylabel{eq:5.6}
\eeq
\beq
\begin{split}
5\psi(q^5)\phi(q^5)^2 & = \psi(q^5)\phi(q^3)^2+4\psi(q^2)\psi(q^3)\phi(q^{15}) \\
& + 8q^4\psi(q)\psi(q^6)\psi(q^{30}) + 4q^5\phi(q^9)^2\psi(q^{45})+8q^7f(q^3,q^{15})^2\psi(q^{45}) \\
& + 8qf(q^{15},q^{30})f(q^3,q^{15})\phi(q^9) + 4q^2f(q^3,q^{15})^2f(q^{15},q^{30}) \\
& - 4q^2\psi(q^6)\psi(q^{15})\phi(q)-4q^7\psi(q)\psi(q^{60})\phi(q^6) \\
& - 4q\psi(q)\psi(q^{12})\phi(q^{30})-4q\phi(q^3)\psi(q^3)\psi(q^{10}), 
\end{split}
\mylabel{eq:5.7}
\eeq
\beq
\begin{split}
15q\psi(q^{15})\phi(q^{15})^2& = -q\psi(q^{15})\phi(q)^2-4\psi(q)\psi(q^6)\phi(q^5) \\  
& + 4\psi(q^2)\psi(q^5)\phi(q^3)+4f(q,q^5)^2f(q^5,q^{10})+4q\psi(q^{15})\phi(q^3)^2  \\
& + 8q\psi(q^2)\psi(q^3)\psi(q^{10})-4q^3\psi(q)\psi(q^{30})\phi(q)-4\psi(q^3)\psi(q^4)\phi(q^{10}) \\
& - 4q^2\psi(q^3)\psi(q^{20})\phi(q^2). 
\end{split}
\mylabel{eq:5.8}
\eeq
Moreover, using some elbow grease one checks that the above is just a generating function form of the following statements.
\beq
3(1,8,8,0,0,0)\left(\frac{M}{3^2}\right)=-W_{1}(M)-W_{2}(M)+W_{3}(M)+W_{4}(M),
\mylabel{eq:5.9}
\eeq
with $M\equiv 1\mymod 8$, $3|M$, 
\beq
5(1,8,8,0,0,0)\left(\frac{M}{5^2}\right)=W_{1}(M)-W_{2}(M)+W_{3}(M)-W_{4}(M),
\mylabel{eq:5.10}
\eeq
with $M\equiv1\mymod 8$, $5|M$, 
\beq
15(1,8,8,0,0,0)\left(\frac{M}{15^2}\right)=-W_{1}(M)+W_{2}(M)+W_{3}(M)-W_{4}(M),
\mylabel{eq:5.11}
\eeq
with $M\equiv 1\mymod 8$, $15|M$.
To state \eqn{5.9}--\eqn{5.11}  in an economical manner we need to develop appropriate notation.
Let $n_{i}$ be some integer represented by $TG_{i}$, such that $\gcd(n_{i},w)=1$ for some $1\leq w$, $w|15$. 
Next, we define $\epsilon(i,w)$ as
\beqs
\epsilon(i,w):= (-n_{i}| w).
\eeqs
We also require that 
\beqs
\epsilon(i,1):= 1.
\eeqs
It is important to realize that this definition does not depend on one's choice of $n_{i}$.
We are now well equipped to combine \eqn{5.2}, \eqn{5.9}--\eqn{5.11} into the single potent statement
\beq
w(1,8,8,0,0,0)\left(\frac{M}{w^2}\right)=\sum_{i=1}^{4}\epsilon(i,w) W_{i}(M), 
\mylabel{eq:5.12}
\eeq
where $w=1,3,5,15$ and $M\equiv1\mymod8$, $w|M$.
As before one can extend \eqn{4.12} by using
$\frac{1}{3}(1,1,1,0,0,0)(M)$ instead of $(1,8,8,0,0,0)(M)$
\beq
w(1,1,1,0,0,0)\left(\frac{M}{w^2}\right)=3\sum_{i=1}^{4}\epsilon(i,w)W_{i}(M), 
\mylabel{eq:5.13}
\eeq
where $w=1,3,5,15$ and $M\equiv 1,2\mymod4$, $w|M$. \\
We forgo the proof.

\bigskip
\section{$S-$genus} \label{sec:sgen}
\medskip

Let $S$ be an odd and square free number and let $S=p_1p_2\ldots p_r$ be the prime factorization of $S$.
In this section we introduce (what we believe to be new) a notion of $S$-genus of ternary forms.
To this end we define an injective map from genera of binary quadratic forms of discriminant $-8S$ to
genera of ternary quadratic forms of discriminant $16S^2$.
According to Theorem 3.15 in \cite{Cox}, there are exactly $2^r$ of these genera of binary quadratic forms 
$BG_1,\ldots ,BG_{2^r}$.
Let $ax^2+bxy+cy^2$ be some quadratic form in $BG_i$, with some $1\le i\le 2^r$. We convert it into a ternary form
$$f(x,y,z): = ax^2+ |b|xy+cy^2 + 2Sz^2.$$ 
Next, we extend $f$ to a genus $TG_i$ that contains $f$.
It can be shown that the map 
\beqs
BG_i\rightarrow TG_i,\quad i=1,2,\ldots\ ,2^{r} 
\eeqs
does not depend on what specific binary form from $BG_i$ we decided to start with. 
We can now define the $S$-genus as a union
\beq
S: = TG_1\cup TG_2\cup\ldots\cup TG_{2^r}. 
\mylabel{eq:6.1}
\eeq

We have an elementary proof that the \(TG_i\)'s are disjoint , using all
the special features of the construction. We wondered if all the structure
of the S-genus were really required for the proof, and the answer is no.
Indeed, Ben Kane kindly sent us a short proof of the following fact: if
\(g_1(x,y)\) and \(g_2(x,y)\) are positive primitive binary quadratic
forms of any common discriminant, and \(N\) is any positive integer, then
\(g_1(x,y) + N z^2\) and \(g_2(x,y) + N z^2\) are in the same genus if and
only if \(g_1(x,y)\) and \(g_2(x,y)\) are in the same genus. He  also
gave counterexamples to the strictly analogous conjecture for
spinor genera.

Let $n_{i}$ be some integer represented by $TG_{i}$ such that $\gcd(n_{i},w)=1 $ for some $1\leq w$, $w|S$.
For any positive divisor $w$ of $S$ we define $\epsilon(i,w)$ as
\beq
\epsilon(i,w):=(-n_{i}|w),
\mylabel{eq:6.2}
\eeq
and for that matter we always take $\epsilon(i,1):= 1$. Again, we remark that this definition does not depend on our choice of $n_{i}$ .
For those with some background in quadratic forms, we comment that for prime divisor $p$ of $S$,
$\epsilon(i,p)=1$ iff the forms of $TG_{i}$ are isotropic over the $p$-adic numbers.
We propose that for $i=1,2,\ldots,2^r$
\beq
M_{i}=\prod_{j=1}^{r}\frac{p_{j}+\epsilon(i,p_{j})}{2}, 
\mylabel{eq:6.3}
\eeq
where
\beq
M_{i}:=\sum_{f\in TG_{i} }\frac{16}{|\mbox{Aut}(f)|}.  
\mylabel{eq:6.4}
\eeq
This seems to generalize Lemma 6.6 on page 152 in \cite{Cassels}.
One way to see why the $S-$genus is such an appealing construct is to consider a mass for the $S$-genus, defined by
\beq
M(S):=\sum_{f\in S}\frac{16}{|\mbox{Aut}(f)|}=M_{1}+\ldots+M_{2^r}.
\mylabel{eq:6.5}
\eeq
Remarkably, \eqn{6.3} together with the orthogonality relation
\beq
\sum_{i=1}^{2^r}\epsilon(i,w)=0,\quad 2\leq w,w|S.
\mylabel{eq:6.6}
\eeq
implies that 
\beq
M(S)=S.   
\mylabel{eq:6.7}
\eeq
Thus, $2^r$ genera conspire to produce a startling simplification.
Perhaps, more important is the fact that all our identities for ternary forms can be stated in laconic fashion as
\beq
(1,1,1,0,0,0)(M)=\sum_{f\in S}\frac{48R_{f}(M)}{|\mbox{Aut}(f)|},
\mylabel{eq:6.8}
\eeq
with $M\equiv 1,2\mymod 4$.
For any $2\leq w,w|S$
\beq
w(1,1,1,0,0,0)\left(\frac{M}{w^2}\right)=3\sum_{i=1}^{2^r}\epsilon(i,w)W_{i}(M),   
\mylabel{eq:6.9}
\eeq
with $M\equiv1,2\mymod4$, $w|M$ and 
\beq
W_{i}(M)=16\sum_{f\in TG_{i}}\frac{R_{f}(M)}{|\mbox{Aut}(f)|}.
\mylabel{eq:6.10}
\eeq
We propose that \eqn{6.8}, \eqn{6.9} hold true for any square free odd $S$.
It follows from \eqn{6.9} that for $M\equiv 1,2\mymod 4$ and $\gcd(S^2,M)=S$,
\beqs
W_{1}(M)=W_{2}(M)=\ldots=W_{2^r}(M).  
\eeqs
We point out that for $M\equiv 1,2\mymod 4$ and $\gcd(S,M)=1$, only a single genus in the $S-$genus is allowed to have     forms that represent $M$. What is remarkable about \eqn{6.8} is that it continues to be true as
$\gcd(M,S)$ increases, and indeed as $M$ becomes divisible by high powers of several $p_{i}$.
The proofs of \eqn{6.3}, \eqn{6.8},\eqn{6.9} will be given elsewhere.

\bigskip
\noindent
\textbf{Acknowledgements}
\medskip

\noindent
We are grateful to Benjamin  Kane for his insights. His Corollary 1.3 was crucial to our investigation.
We would like to thank Manjul Bhargava, Frank Garvan, Jonathan Hanke, Byungchan Kim, Michael Somos  
for their kind interest and helpful discussions.
\bigskip
\bigskip

\bibliographystyle{amsplain}

\end{document}